\documentclass[a4paper,11pt]{article}
\usepackage[utf8]{inputenc}
\usepackage[T1]{fontenc}
\usepackage[english]{babel}

\usepackage[intlimits]{amsmath}
\usepackage{amssymb, amsfonts}
\usepackage[amsmath, thmmarks]{ntheorem}

\usepackage{MnSymbol}
\usepackage{mathbbol}
\usepackage{bbm}
\usepackage{mathrsfs}
\usepackage[all]{xy}

\usepackage{graphicx}
\numberwithin{equation}{section}

\newtheorem{theoremcounter}{theoremcounter}[section]

\newtheorem{corollary}[theoremcounter]{Corollary}
\newtheorem{definition}[theoremcounter]{Definition}

\newtheorem{lemma}[theoremcounter]{Lemma}
\newtheorem{remark}[theoremcounter]{Remark}
\newtheorem{theorem}[theoremcounter]{Theorem}

\theorembodyfont{\normalfont}
\theoremsymbol{\ensuremath{\Box}}
\newtheorem*{proof}{Proof}

\newcommand{\NN}{\ensuremath{\mathbb{N}}}
\newcommand{\ZZ}{\ensuremath{\mathbb{Z}}}

\newcommand{\CC}{\ensuremath{\mathbb{C}}}

\newcommand{\Hom}{\ensuremath{\textup{Hom}}}
\newcommand{\id}{\ensuremath{\textup{id}}}

\newcommand{\mat}[2]{\ensuremath{\textup{M}_{#1}(#2)}}
\newcommand{\Cmat}[1]{\ensuremath{\textup{M}_{#1}(\CC)}}
\newcommand{\Cstar}{\ensuremath{C^*}}

\newcommand{\cont}{\ensuremath{\mathcal{C}}}
\newcommand{\cspan}{\ensuremath{\overline{\textup{span}}}}
\newcommand{\Rep}{\ensuremath{\textup{Rep}}}
\newcommand{\Repirreven}{\ensuremath{\textup{Rep}_\textup{even}^\textup{irr}}}
\newcommand{\Repirrodd}{\ensuremath{\textup{Rep}_\textup{odd}^\textup{irr}}}
\newcommand{\Repirrevenone}{\ensuremath{\textup{Rep}_\textup{even}^{\textup{irr}, (1)}}}
\newcommand{\Repirreventwo}{\ensuremath{\textup{Rep}_\textup{even}^{\textup{irr}, (2)}}}
\newcommand{\Repirroddone}{\ensuremath{\textup{Rep}_\textup{odd}^{\textup{irr}, (1)}}}
\newcommand{\Repirroddtwo}{\ensuremath{\textup{Rep}_\textup{odd}^{\textup{irr}, (2)}}}

\newcommand{\Seven}{\ensuremath{S_{\textup{even}}}}
\newcommand{\Sodd}{\ensuremath{S_{\textup{odd}}}}
\newcommand{\Sevenone}{\ensuremath{S_{\textup{even}}^{(1)}}}
\newcommand{\Seventwo}{\ensuremath{S_{\textup{even}}^{(2)}}}
\newcommand{\Soddone}{\ensuremath{S_{\textup{odd}}^{(1)}}}
\newcommand{\Soddtwo}{\ensuremath{S_{\textup{odd}}^{(2)}}}

\newcommand{\PART}{\ensuremath{\textup{Part}}}
\newcommand{\NC}{\ensuremath{\textup{NC}}}

\newcommand{\Aon}{\ensuremath{\textup{A}_\textup{o}(n)}}
\newcommand{\Abn}{\ensuremath{\textup{A}_\textup{b}(n)}}
\newcommand{\Abprimen}{\ensuremath{\textup{A}_\textup{b'}(n)}}
\newcommand{\Asn}{\ensuremath{\textup{A}_\textup{s}(n)}}
\newcommand{\Asprimen}{\ensuremath{\textup{A}_\textup{s'}(n)}}
\newcommand{\Ahn}{\ensuremath{\textup{A}_\textup{h}(n)}}
\newcommand{\Astarn}{\ensuremath{\textup{A}_*(n)}}
\newcommand{\Aun}{\ensuremath{\textup{A}_\textup{u}(n)}}
\newcommand{\Akn}{\ensuremath{\textup{A}_\textup{k}(n)}}
\newcommand{\Acn}{\ensuremath{\textup{A}_\textup{c}(n)}}
\newcommand{\Apn}{\ensuremath{\textup{A}_\textup{p}(n)}}

\newcommand{\mon}[1]{\ensuremath{{\textup{mon}(#1)}}}

\begin{document}

% \subjclass[2000]{Primary 46L65.}
% \keywords{Quantum groups, Fusion rules.}

\begin{center}
{\LARGE \bf {Isomorphisms and Fusion Rules of Orthogonal Free Quantum Groups and their Free Complexifications}}
\end{center}
\bigskip
\begin{center}
{\sc {by Sven Raum$^{(1,2)}$}
\setcounter{footnote}{1}
\footnotetext{Research partially supported by Marie Curie Research Training Network Non-Commutative Geometry MRTN-CT-2006-031962 and by K.U.Leuven BOF research grant OT/08/032}
\setcounter{footnote}{2}
\footnotetext{Department of Mathematics;    K.U.Leuven; Celestijnenlaan 200B; B--3001 Leuven (Belgium).    \\ E-mail: sven.raum@wis.kuleuven.be}}
\end{center}

\subsection*{Abstract}
We show that all orthogonal free quantum groups are isomorphic to variants of the free orthogonal Wang algebra, the hyperoctahedral quantum group or the quantum permutation group. We also obtain a description of their free complexification. In particular we complete the calculation of fusion rules of all orthogonal free quantum groups and their free complexifications.
\section{Introduction}\label{section_introduction}

One problem in the theory of compact quantum groups is to find examples whose invariants can be calculated. The fusion rules of a compact quantum group are one of these invariants. Fusion rules give a complete description of equivalence classes of irreducible corepresentations and a decomposition of the tensor product of two of them into irreducible corepresentations. One approach to this problem is given by 'free quantum groups' as defined in \cite{banicaspeicher09}. These are orthogonal quantum groups, i.e. subgroups of the free orthogonal Wang algebra, whose intertwiners can be described by non-crossing partitions.\\
Given natural numbers $k$ and $l$ the set $\PART(k,l)$ denotes the set of all partitions on two rows with $k$ and $l$ points, respectively. That is, an element $P \in \PART(k,l)$ is a partition of the disjoint union $\{1,...,k\} \sqcup \{1,...,l\}$. Alternatively it can be described by a diagram
\begin{equation*}
\left\{\begin{xy}
	\POS(0,5) 	*{\cdot}
	\POS(4,5) 	*{\cdot}
	\POS(8,5) 	*{...}
	\POS(12,5)		*{\cdot}
	
	\POS(6,0)		*{P}
	
	\POS(0,-5) 	*{\cdot}
	\POS(4,-5) 	*{\cdot}
	\POS(8,-5) 	*{...}
	\POS(12,-5)	*{\cdot}
\end{xy}\right\}
\end{equation*}
connecting the $k$ points in the upper row and the $l$ points in the lower row according to the partition of $\{1,...,k\} \sqcup \{1,...,l\}$. $P$ is called non-crossing if it can be represented by a diagram with no lines crossing. The set of all non-crossing partitions on $k$ and $l$ points is denoted by $\NC(k,l)$.\\
Let $n, k, l \in \NN$ and let $(e_i)$ be the standard basis of $\CC^n$. Let $i = (i_1,...,i_k) \in \{1,...,n\}^k$ and $j = (j_1,...,j_l) \in \{1,...,n\}^l$ be multi indices and $P \in \PART(k,l)$. We set $P(i,j) = 1$ if and only if the diagram $P$ joins only equal numbers after writing the entries of $i$ in the upper row of the above diagram and those of $j$ in the lower row. If $P$ connects different numbers set $P(i,j) = 0$.\\
Using this notation, a partition $P \in \PART(k,l)$ defines a linear map $T_P$ from $(\CC^n)^{\otimes k}$ to $(\CC^n)^{\otimes l}$ by \[T_P(e_{i_1} \otimes ... \otimes e_{i_k}) = \sum_{j_1, ..., j_l}{P(i_1,...,i_k;j_1,...,j_l) \cdot e_{j_1} \otimes ... \otimes e_{j_l}}.\] A subspace of $\Hom((\CC^n)^{\otimes k}, (\CC^n)^{\otimes l})$ is by definition spanned by partitions if it is linearly generated by a family $(T_P)$ where $P$ runs through some subset of $\PART(k,l)$.\\
In \cite{wang95} the free unitary Wang algebra
\[A_u(n) := \Cstar(u_{ij}, \; 1 \leq i,j \leq n | (u_{ij})_{ij}, (u_{ij}^*)_{ij} \text{ are unitary})\]
and the free orthogonal Wang algebra
\[A_o(n) := \Cstar(u_{ij}, \; 1 \leq i,j \leq n | (u_{ij})_{ij} = (u_{ij}^*)_{ij} \text{ is unitary})\]
were introduced. Moreover in \cite{wang98} the quantum permutation group 
\[A_s(n) := \Cstar\left(u_{ij}, \; 1 \leq i,j \leq n \middle|
\begin{array}{l} (u_{ij}) = (u_{ij}^*) \text{ is unitary and } u_{ij} \text{ are} \\
\text{partial isometries summing up to one} \\
\text{in every row and every column} \end{array}\right)\]
was defined. Note that ``are partial isometries'' can be replaced by ``are projections''. The three last named algebras are compact matrix quantum groups in the sense of Woronowicz \cite{woronowicz91}. \\
The following class of quantum groups will be of interest in this paper.

\begin{definition}\label{definition_free_quantum_groups}
Let $(A,U)$ be a compact matrix quantum group. Then it is called free if
\begin{itemize}
  \item The morphism $(\Aun, U_\textup{u}) \rightarrow (\Asn, U_\textup{s})$ mapping the entries of $U_\textup{u}$ to those of $U_\textup{s}$ factorizes through $(A,U)$.
  \item The intertwiner spaces $\Hom(U^{i_1} \boxtimes \dotsm \boxtimes U^{i_k}, U^{j_1}\boxtimes \dotsm \boxtimes U^{j_l})$, $i_\alpha, j_\beta \in \{1, \overline{\phantom{U}}\}$ are spanned by partitions, where $\overline{U} = (u_{ij}^*)$ is the conjugate corepresentation of $U$ and $\boxtimes$ denotes the tensor product of corepresentations.
\end{itemize}
If the first condition is strengthened by requiring that the morphism $(\Aon, U_\textup{o}) \rightarrow (\Asn, U_\textup{s})$ factors through $(A,U)$, then $A$ it is called orthogonal free.
\end{definition}

In \cite{banicaspeicher09} the following classification was achieved.

\begin{theorem}
  There are exactly six orthogonal free quantum groups. Namely
  \begin{enumerate}
    \item The free orthogonal Wang algebra.
    \item The quantum permutation group.
    \item The hyperoctahedral quantum group
    	\[\Ahn := \Cstar\left(u_{ij}, \; 1 \leq i,j \leq n \middle|
    	\begin{array}{l} (u_{ij}) = (u_{ij}^*) \text{ is unitary and } \\
    	u_{ij} \text{ are partial isometries} \end{array}\right).\]
    \item The bistochastic quantum group 
    	\[\Abn := \Cstar\left(u_{ij}, \; 1 \leq i,j \leq n \middle|
 			\begin{array}{l} (u_{ij}) = (u_{ij}^*) \text{ is unitary and } \\
 			u_{ij} \text{ sum up to one} \\
 			\text{in every row and every column} \end{array}\right).\]
    \item The symmetrized bistochastic quantum group
    	\[\Abprimen := \Cstar\left(u_{ij}, \; 1 \leq i,j \leq n \middle|
 			\begin{array}{l} (u_{ij}) = (u_{ij}^*) \text{ is unitary and } \\
 			u_{ij} \text{ sum up to the same element} \\
 			\text{in every row and every column} \end{array}\right).\]
    \item The symmetrized quantum permutation group
    	\[\Asprimen := \Cstar\left(u_{ij}, \; 1 \leq i,j \leq n \middle|
 			\begin{array}{l} (u_{ij}) = (u_{ij}^*) \text{ is unitary and } \\
 			u_{ij} \text{ are partial isometries} \\
 			\text{summing up to the same element} \\
 			\text{in every row and every column} \end{array}\right).\]
  \end{enumerate}
\end{theorem}

The fusion rules of (1) were calculated in \cite{banica96}, those of (2) in \cite{banica99} and those of (3) in \cite{banicavergnioux09}. We show that the remaining examples are slight modifications of $\Aon$ and $\Asn$. In particular we can derive their fusion rules and find that $\Abprimen$ and $\Asprimen$ are counterexamples to a conjecture by Banica and Vergnioux given in \cite{banicavergnioux09}.\\
In \cite{banica08} the free complexification of orthogonal free quantum groups was considered. If $(A,U)$ is a orthogonal free quantum group, then its free complexification $(\widetilde{A}, \widetilde{U})$ is by definition the sub-C*-algebra of the free product $A * \cont(S^1)$ generated by the entries of $\widetilde{U} := U \cdot \id_{\textup{S}^1} = (u_{ij} \cdot \id_{\textup{S}^1})$. Here $\id_{\textup{S}^1}$ denotes the canonical generator of $\cont(\textup{S}^1)$. As Banica shows in \cite{banica08} the intertwiners between tensor products of the fundamental corepresentation and its conjugate can be described by the intertwiners of the orthogonal free quantum group it comes from. With additional requirements we can calculate the fusion rules of the free complexification from the fusion rules of the original orthogonal free quantum group. These additional requirements are fulfilled by $\Aon$ and $\Ahn$, which gives the fusion rules of $\Akn = \widetilde{\Ahn}$. Those of $\Aun = \widetilde{\Aon}$ are known from \cite{banica97}. \\
From \cite{banica08} we know that $\widetilde{\Abn} = \widetilde{\Abprimen}$ and $\widetilde{\Asn} = \widetilde{\Asprimen}$. We denote $\widetilde{\Abn} =: \Acn$ and $\widetilde{\Asn} =: \Apn$. They can be decomposed and described in terms of $\Aon$ and $\Asn$ again. \\
\textbf{Acknowledgment:} I want to thank both Thomas Timmermann for suggesting to work on fusion rules of free quantum groups and Stefaan Vaes for helpful discussions about this article, especially on the last section. Moreover, I want to thank the referee for helpful comments.

\section{Preliminaries}

We will mainly work with compact matrix quantum groups as defined by Worono-wicz in \cite{woronowicz91}. If $A$ is a *-algebra and $U \in \mat{n}{A}$ we denote by $\overline{U}$ the matrix whose entries are conjugated, i.e. $\overline{U}_{ij} = (U_{ij})^*$. \\
A pair $(A,U)$ of a C*-algebra $A$ and a unitary $U \in \mat{n}{A}$ is called a compact matrix quantum group if
\begin{itemize}
 \item $A$ is generated by the entries of $U$,
 \item there is a *-homomorphism $\Delta:A \rightarrow A \otimes_\textup{min} A$ mapping $u_{ij}$ to $\sum_k u_{ik} \otimes u_{kj}$,
 \item the matrix $\overline{U}$ is invertible.
\end{itemize}

A morphism of compact matrix quantum groups $(A,U) \stackrel{\phi}{\rightarrow} (B,V)$ is a *-homo-morphism $A \rightarrow B$ such that $\phi(u_{ij}) = v_{ij}$ where $U$ and $V$ must have the same size. There is at most one morphism from one quantum group to another. If there is a morphism $(A,U) \rightarrow (B,V)$ then we say that $(B,V)$ is a quantum subgroup of $(A,U)$. \\
Every compact matrix quantum group is also a compact quantum group, i.e. a C*-algebra $A$ with a *-homomorphism $\Delta:A \rightarrow A \otimes_\textup{min} A$ such that
\begin{itemize}
 \item $(\Delta \otimes \id) \circ \Delta = (\id \otimes \Delta) \circ \Delta$,
 \item $\cspan(A \otimes 1)\Delta(A) = \cspan(1 \otimes A)\Delta(A) = A \otimes A$.
\end{itemize}
A morphism of compact quantum groups $(A,\Delta_A) \stackrel{\phi}{\rightarrow} (B,\Delta_B)$ is a *-homomorphism from $A$ to $B$ such that $\Delta _B\circ \phi = (\phi \otimes \phi) \circ \Delta_A$.
Every morphism of compact matrix quantum groups is also a morphism of compact quantum groups.\\
We will also refer to a quantum group $(A,U)$ or $(A,\Delta)$ as $A$. \\
If $(A,\Delta_A)$ and $(B,\Delta_B)$ are quantum groups, then we denote by $(A,\Delta_A) \otimes (B,\Delta_B)$ the direct sum of quantum groups and by $(A,\Delta_A) * (B,\Delta_B)$ their free product. We will also write $A \otimes B$ and $A * B$. \\
A unitary corepresentation matrix of $(A,\Delta)$ is a unitary matrix $V \in \mat{m}{A}$ such that $\Delta(v_{ij}) = \sum_k v_{ik} \otimes v_{kj}$. In particular a one dimensional corepresentation matrix is just a unitary group-like element of $A$.
\section{Free fusion rings}\label{section_fusion_rings}				

In this section we will introduce free fusion rings and prove that they are free unital rings. \\

We will use the following notation for words in free monoids. Let $M = \mon{S}$ be a free monoid over a set $S$. If $w \in M$ is a word of length $k$, then we write $w_i$ for the $i$-th letter of $w$, $1 \leq i \leq k$. Hence $w = w_1 w_2 w_3 \dotsc w_{k-1} w_k$.

\begin{definition}\label{definition_free_fusion_monoid}
A free fusion monoid is a free monoid $M = \mon{S}$ over a set $S$ with a fusion $\cdot \,:S \times S \rightarrow S \cup \{\emptyset\}$ and a conjugation $\overline{\phantom{u}}:S \longrightarrow S$. They must satisfy the following conditions.
\begin{enumerate}
  \item The fusion $\cdot$ is associative, where we make the convention that $s \cdot s'$ is the empty set if one of $s,s'$ is the empty set.
  \item The conjugation is involutive, i.e. $\overline{\overline{s}} = s$ for all $s \in S$.
  \item Fusion and conjugation are compatible in the following sense. For all $s_1, s_2, s_3 \in S$ we have
    \[s_1 \cdot s_2 = \overline{s_3} \Leftrightarrow s_2 \cdot s_3 = \overline{s_1}\]
\end{enumerate}
A set $S$ equipped with fusion and conjugation is called a fusion set. \\
The fusion and conjugation of $S$ induce a fusion and a conjugation on $M$ via
\begin{itemize}
  \item $w \cdot w' = w_1 \dotsc w_{k-1}(w_k \cdot w'_1)w'_2 \dotsc w'_l$ where this fusion is the empty set by convention if $w_k \cdot w'_1 = \emptyset$.
  \item $\overline{w} = \overline{w_k} \dotsc \overline{w_1}$
\end{itemize}
\end{definition}

If $M = \mon{S}$ is a free fusion monoid, we can turn $\ZZ M$ into an associative ring by
\[a_w \cdot a_{w'} = \sum_{\substack{w = xy \\ w' = \overline{y}z}}{(a_{xz } + a_{x \cdot z})}.\]
Here $w$, $w'$ are words in $M$, $a_w$ and $a_{w'}$ are the corresponding elements in $\ZZ M$, $xy$, $\overline{y}z$ and $xz$ denote the concatenation of words and the second term in the sum is by convention always ignored if the fusion $x \cdot z$ is empty. Actually condition (3) of the previous definition is a necessary condition for making $\ZZ M$ associative, as it can be seen by considering $(a_{s_1} \cdot a_{s_2}) \cdot a_{s_3} = a_{s_1} \cdot ( a_{s_2} \cdot a_{s_3})$ for $s_1, s_2, s_3 \in S$. A *-ring isomorphic to $\ZZ M$ for some fusion monoid $M$ is called a free fusion ring. \\
From the point of view of rings, free fusion rings are very easy. Actually they are free. The proof of the following lemma was already given in \cite{banicavergnioux09} in some special cases.

\begin{lemma}\label{lemma_free_fusion_rings_are_free}
A free fusion ring over a fusion set $S$ is the free unital ring over $a_s$, $s \in S$.
\begin{proof}
  Let $\ZZ M$ be the fusion ring over a fusion set $S$. It suffices to show that $\ZZ M$ is a free $\ZZ$-module with the basis $a_{s_1} \dotsm a_{s_k}$ with  $k \in \NN$ and $s_1, \dotsc ,s_k \in S$. So it suffices to express the elements of the $\ZZ$-basis $a_w, w \in M$ as $\ZZ$-linear combinations of the elements $a_{s_1} \dotsm a_{s_k}$ with  $k \in \NN$ and $s_1, \dotsc ,s_k \in S$ and to show that $\{a_{s_1} \dotsm a_{s_k} |  k \in \NN, s_1, \dotsc ,s_k \in S\}$ is $\ZZ$-linearly independent. \\
  There are coefficients $C_{s_1 \dotsc s_k}^w \in \ZZ$ such that $a_{s_1} \dotsm a_{s_k} = a_{s_1 \dotsc s_k} + \sum_{|w| < k}{C_{s_1 \dotsc s_k}^w a_w}$, where $|w|$ is the length of the word $w \in M$. This shows that  $\{a_{s_1} \dotsm a_{s_k} |  k \in \NN, s_1, \dotsc ,s_k \in S\}$ is linearly independent. Moreover, by induction on $k$ there are coefficients $D_{s_1 \dotsc s_l}^w \in \ZZ$ such that $a_{s_1 \dotsc s_k} = a_{s_1} \dotsm a_{s_k} + \sum_{|w| < k}{D_{s_1 \dotsc s_k}^w a_{w_1} \dotsm a_{w_{|w|}}}$. This shows that all $a_w, w \in M$ are linear combinations of $a_{s_1} \dotsm a_{s_k}$ with  $k \in \NN$ and $s_1, \dotsc,s_k \in S$.
\end{proof}
\end{lemma}

\begin{remark}
 Free fusion rings can be used to describe fusion rules very shortly and there is hope to use free fusion rings as a starting point for proofs of several properties of quantum groups. See section 10 of \cite{banicavergnioux09} for a comment on these possibilities. However in order to justify the concept of free fusion rings intrinsically it would be good to answer the following question affirmatively. Is every fusion ring of a compact quantum group that is free as a unital ring a free fusion ring?
\end{remark}
\section{Some isomorphisms of combinatorial quantum groups}\label{section_isomorphisms}

In this section we will consider combinatorial quantum groups $\Astarn$ for $* \in \{b,b',s',c,p\}$. They are free products  or direct sums of known quantum groups. For $* \in \{b',s',c,p\}$ it turns out that their fusion rings are not free.

\begin{theorem}\label{theorem_isomorphisms}
We have the following isomorphisms of compact quantum groups (not necessarily preserving the fundamental corepresentation).
\begin{enumerate}
  \item $\Abn$ is isomorphic to $\textup{A}_\textup{o}(n - 1)$.
  \item $\Asprimen$ is isomorphic to the direct sum $\Asn \otimes C^*(\ZZ/{2\ZZ})$.
  \item $\Abprimen$ is isomorphic to the free product $\Abn * C^*(\ZZ/{2\ZZ})$.
  \item $\Apn$ is isomorphic to the free product $\Asn * \cont(\textup{S}^1)$.
  \item $\Acn$ is isomorphic to the free product $\Abn * \cont(\textup{S}^1)$.
\end{enumerate}
\end{theorem}

\begin{remark}
 Note that in the case $n \leq 3$ we have the isomorphisms  $\Asn \cong \cont{(\textup{S}_\textup{n})}$ and $\textup{A}_\textup{o}(1) \cong \cont(\{-1,1\})$. So the given descriptions can be further simplified.
\end{remark}

Theorem \ref{theorem_isomorphisms}(1) is proven by the following remark. Let $U \in \mat{n}{A}$ be an orthogonal matrix, i.e. $\overline{U} = U$ unitary, where $A$ is any unital $\Cstar$-algebra. Then $U$ is bistochastic if and only if the vector $(1,1, \dotsc,1)^\textup{t}$ is a right eigenvector and $(1,1, \dotsc,1)$ is a left eigenvector of $U$. If $T \in \Cmat{n}$ denotes any orthogonal matrix such that $T(1,0, \dotsc,0)^\textup{t} = (1/{\sqrt{n}},\dotsc,1/{\sqrt{n}})^\textup{t}$ , then an orthogonal matrix $U$ is bistochastic if and only if $T^\textup{t}UT$ is of block form with $1$ in the upper left corner and an orthogonal $(n - 1) \times (n - 1)$ matrix in the lower right corner. \\
The key observation for the rest of \ref{theorem_isomorphisms} is the following lemma.

\begin{lemma}
Let $* \in \{b', s', c, p\}$. The fundamental corepresentation of $\Astarn$ contains a one dimensional non-trivial corepresentation $U_z$ which fulfils $U_z \boxtimes \overline{U_z} \simeq 1$. If $* \in \{b', s'\}$ then $U_z \simeq \overline{U_z}$. 
\begin{proof}
  Consider $* = b', s'$ first. The element $z = \sum_i{u_{ij}}$ is easily seen to be a unitary group-like element, so it corresponds to a one dimensional unitary corepresentation of $\Astarn$. Consider the group $\textup{S}_n \oplus \ZZ/{2\ZZ} \subset U_n$ as permutation matrices with entries $+1$ and $-1$. Let $U_{\textup{S}_n \oplus \ZZ/{2\ZZ}}$ be the canonical fundamental corepresentation of $\cont(\textup{S}_n \oplus \ZZ/{2\ZZ})$. Then the image of $z$ under the map $(\Astarn, U_*) \rightarrow (\cont(\textup{S}_n \oplus \ZZ/{2\ZZ}), U_{\textup{S}_n \oplus \ZZ/{2\ZZ}})$ is $-1$, so $z$ is non-trivial. \\
  For $* = p,c$ consider $z := \id_{\textup{S}^1}$ as coming from the copy of $\cont(\textup{S}^1)$. This copy is contained in $\Astarn$, since the trivial corepresentation is contained in the fundamental corepresentation of $\Abn$ and $\Asn$. \\
  Using the relations of $\Astarn$ we can check the rest of the claim by simple calculations.
\end{proof}
\end{lemma}

\begin{remark}
The last lemma shows, that the fusion rules of neither of the quantum groups $\Astarn$ for $\* \in \{b', s', c, p\}$ can be described by a free fusion ring. Actually in a free fusion ring any element $a \neq 1$ satisfies $a \cdot a^* \neq 1$. This gives two counterexamples to the conjecture that for $n \geq 4$ the fusion rules of all free orthogonal quantum groups  can be described by a free fusion ring, which was stated in \cite{banicavergnioux09}.
\end{remark}

\begin{remark}
The fundamental corepresentation of any matrix quantum group that has $(\Asn, U_\textup{s})$ as a sub quantum group cannot be the sum of more than two irreducible corepresentations. In particular the last lemma already gives a decomposition $U \simeq U_z \boxplus V$ with $U_z$ non-trivial and one dimensional and $V$ irreducible, where $U$ is the fundamental corepresentation of $\Astarn$.
\end{remark}

\begin{proof}[Proof of Theorem \ref{theorem_isomorphisms}]
 The isomorphism of (2) is given by $\Asn \otimes C^*(\ZZ/{2\ZZ}) \rightarrow \Asprimen: \: u^\textup{s}_{ij} \otimes 1 \mapsto u^{\textup{s}'}_{ij} \cdot z, \, 1 \otimes u_{\overline{1}} \mapsto z$. This map exists since $z$ is central in $\Asprimen$ as an easy calculation shows. The inverse map is given by \[\Asprimen \rightarrow \Asn \otimes C^*(\ZZ/{2\ZZ}): \: u^{\textup{s'}}_{ij} \rightarrow u^{\textup{s}}_{ij} \otimes u_{\overline{1}}.\]
 In order to prove (3) we use again an orthogonal matrix $T \in \Cmat{n}$ such that $T(1,0,...,0)^\textup{t} = (1/{\sqrt{n}},...,1/{\sqrt{n}})^\textup{t}$. Then a matrix $U \in \mat{n}{A}$ for some $\Cstar$-algebra $A$ satisfies the relations of $U_{\textup{b'}}$ if and only if $T^{\textup{t}}UT$ is a block matrix with a self-adjoint unitary in the upper left corner and an orthogonal $(n - 1) \times (n - 1)$ matrix in the lower right corner. This proves $\Abprimen \cong \textup{A}_\textup{o}(n - 1) * C^*(\ZZ/{2\ZZ}) \cong \Abn * C^*(\ZZ/{2\ZZ})$. \\
 The isomorphism of (4) is given by  \[\Asn * \cont(\textup{S}^1) \rightarrow \Apn: \: u^{\textup{s}}_{ij} \mapsto u^\textup{p}_{ij} \cdot z^*, \, \id_{\textup{S}^1} \mapsto z.\]
 The isomorphism of (5) is given by \[\Abn * \cont(\textup{S}^1) \rightarrow \Acn: \: u^{\textup{b}}_{ij} \mapsto u^\textup{c}_{ij} \cdot z^*, \, \id_{\textup{S}^1} \mapsto z.\]
 All the isomorphisms respect the comultiplication, since $z$ is group-like. Hence, they are isomorphisms of quantum groups.
\end{proof}

\section{Fusion rules for free products and the quantum group $\Akn$}

In this section we describe the fusion rules of the free complexification $\Akn \cong \widetilde{\Ahn}$. Instead of referring to $\Akn$ explicitly, we will work in a more general setting and deduce its fusion rules as a corollary. Roughly the main statement of this section is given by the following theorem. See theorem \ref{theorem_fusion_rules_of_complexifications} for a precise statement.

\begin{theorem}\label{theorem_main_statement_simplified}
  Let $(A,U)$ be an orthogonal compact matrix quantum group, i.e. $\overline{U} = U$, such that its fusion rules are free. Assume further that $1 \notin U^{\boxtimes 2k + 1}$ for any $k \in \NN$. Then the fusion rules of $(\widetilde{A}, \widetilde{U})$ are free and can be described in terms of the fusion rules of $(A,U)$.
\end{theorem}

The following theorem is due to Wang \cite{wang98}.
\begin{theorem}\label{theorem_irr_rep_of_free_products}
  Let $(A,\Delta_A)$ and $(B, \Delta_B)$ be compact quantum groups. Let $(U^\alpha)_{\alpha \in \mathscr{A}}$ and $(U^\beta)_ {\beta \in \mathscr{B}}$ be complete sets of representatives of irreducible corepresentations of $A$ and $B$, respectively. Then the corepresentations $(W^{\gamma_1} \boxtimes \dotsm \boxtimes W^{\gamma_n})$ with $n \in \NN$, all $W^{\gamma_i}$ in $\{U^\alpha \, | \, \alpha \in \mathscr{A}\}$ and $\{U^\beta \, | \, \beta \in \mathscr{B}\}$ and neighbours not from the same set, form a complete set of irreducible representations of the free product $(A, \Delta_A) * (B, \Delta_B)$. 
\end{theorem}

The following observation will be useful when studying the fusion rules of a free complexification.

\begin{remark}\label{remark_fusion_rules_in_free_products}
  Let $A * B$ be a free product of compact quantum groups with irreducible corepresentations $W^{\gamma_1} \boxtimes \dotsm \boxtimes W^{\gamma_n}$ and $W^{\delta_1} \boxtimes \dotsm \boxtimes W^{\delta_m}$ as in the last theorem. Then
  \begin{enumerate}
    \item If $W^{\gamma_n}$ and $W^{\delta_1}$ are not corepresentations of the same factor of the free product, then $W^{\gamma_1} \boxtimes \dotsm \boxtimes W^{\gamma_n} \boxtimes W^{\delta_1} \boxtimes \dotsm \boxtimes W^{\delta_m}$ is an irreducible corepresentation of $A * B$.
    \item If $W^{\gamma_n}$ and $W^{\delta_1}$ are corepresentations of the same factor and $W^{\gamma_n} \boxtimes W^{\delta_1} = \sum_{i = 1}^k W^{\epsilon_i} + \delta_{\overline{W^{\gamma_n}}, W^{\delta_1}} \cdot 1$ is the decomposition into irreducible corepresentations, then
    \begin{multline*}
      W^{\gamma_1} \boxtimes \dotsm \boxtimes W^{\gamma_n} \boxtimes W^{\delta_1} \boxtimes \dotsm \boxtimes W^{\delta_m} \\
      = \sum_{i = 1}^k{(W^{\gamma_1} \boxtimes \dotsm \boxtimes W^{\gamma_{n-1}} \boxtimes W^{\epsilon_i} \boxtimes W^{\delta_2} \boxtimes \dotsm \boxtimes W^{\delta_m})} \\
      + \delta_{\overline{W^{\gamma_n}}, W^{\delta_1}} \cdot W^{\gamma_1} \boxtimes \dotsm \boxtimes W^{\gamma_{n-1}} \boxtimes W^{\delta_2} \boxtimes \dotsm \boxtimes W^{\delta_m}
    \end{multline*}
    and the first $k$ summands of this decomposition are irreducible.
  \end{enumerate}
\end{remark}

For the rest of this section fix an orthogonal compact matrix quantum group $(A,U)$ such that its fusion rules are described by a free fusion ring over the fusion set $S$. Assume further that $1 \notin U^{\boxtimes 2k + 1}$ for any $k \in \NN$. \\
Note that the fusion ring of $\widetilde{A}$ is the fusion subring of $\Rep(A * \cont(\textup{S}^1))$ that is generated by $U \boxtimes z$, where $z$ denotes the identity on the circle.\\
We will construct the free complexification $\widetilde{S}$ of $S$ and prove that the fusion rules of $(\widetilde{A}, \widetilde{U})$ are described by $\widetilde{S}$. We begin by constructing $\widetilde{S}$. \\
Let $\Repirreven$ (respectively $\Repirrodd$) be the set of classes of irreducible corepresentations of $A$ that appear as subrepresentations of an even (respectively odd) tensor power of $U$. We have $\Repirreven \cap \Repirrodd = \emptyset$ due to Frobenius duality and the requirement $1 \notin U^{2k + 1}$ for all $k \in \NN$. Let $\Seven \subset S$ (resp. $\Sodd \subset S$) be the set of elements corresponding to corepresentations from $\Repirreven$ (resp. $\Repirrodd$). The set $\widetilde{S}$ is then by definition the disjoint union $\Seven \sqcup \Seven \sqcup \Sodd \sqcup \Sodd$. Denote the first copy of $\Seven$ (resp. $\Sodd$) by $\Sevenone$ (resp. $\Soddone$) and the second one by $\Seventwo$ (resp. $\Soddtwo$). \\

What follows is motivated by the following point of view:
\begin{remark}\label{remark_elements_of_S_tilde}
  We consider element of $\Sevenone$ as a plain copy of those in $\Seven$. The elements of $\Seventwo$ are of the form $z^*\cdot s \cdot z$ for some $s \in \Seven$. Similarly we consider elements of $\Soddone$ as $s \cdot z$ and elements of $\Soddtwo$ as $z^* \cdot s$ for $s \in \Sodd$.
\end{remark}

Define a conjugation on $\widetilde{S}$ by the conjugation on $S$ leaving $\Sevenone$ and $\Seventwo$ globally invariant and exchanging $\Soddone$ and $\Soddtwo$. Note that $\overline{\Seven} = \Seven$ and $\overline{\Sodd} = \Sodd$, i.e. the conjugation on $\widetilde{S}$ is well defined. A fusion on $\widetilde{S}$ can be defined according to the following table.
\begin{center}
\begin{tabular}{c||c|c|c|c}
  $\cdot$	&\Sevenone	&\Seventwo	&\Soddone	&\Soddtwo	\\
  \hline
  \Sevenone	&$\Sevenone \cup \{\emptyset\}$	&$\emptyset$	&$\Soddone \cup \{\emptyset\}$	&$\emptyset$	\\
  \hline
  \Seventwo	&$\emptyset$	&$\Seventwo \cup \{\emptyset\}$	&$\emptyset$	&$\Soddtwo \cup \{\emptyset\}$	\\
  \hline
  \Soddone	&$\emptyset$	&$\Soddtwo \cup \{\emptyset\}$	&$\emptyset$	&$\Sevenone \cup \{\emptyset\}$	\\
  \hline
  \Soddtwo	&$\Soddtwo \cup \{\emptyset\}$	&$\emptyset$	&$\Seventwo \cup \{\emptyset\}$	&$\emptyset$
\end{tabular}
\end{center}
The row gives the element which is fused from the right with an element coming from the set indicated by the column. The fusion is empty if this is indicated by the table and is otherwise the usual fusion of two elements of $S$ lying in the part of $\widetilde{S}$ indicated by the table. Note that this definition makes sense, since $\Seven \cdot \Seven, \Sodd \cdot \Sodd \subset \Seven \cup \{\emptyset\}$ and $\Seven \cdot \Sodd, \Sodd \cdot \Seven \subset \Sodd \cup \{\emptyset\}$. It is easy to see that $\widetilde{S}$ with this structure is a fusion set. \\
Now we can state a precise version of \ref{theorem_main_statement_simplified}.

\begin{theorem}\label{theorem_fusion_rules_of_complexifications}
  Let $(A,U)$ be an orthogonal compact matrix quantum group such that its fusion rules are described by a free fusion ring over the fusion set $S$. Assume further that $1 \notin U^{\boxtimes 2k + 1}$ for any $k \in \NN$. Then the fusion rules of $(\widetilde{A}, \widetilde{U})$ are given by the free complexification $\widetilde{S}$ of $S$.
\end{theorem}

We construct a complete set of corepresentations of $\widetilde{A}$. In order to do so we associate an irreducible corepresentations of $(\widetilde{A}, \widetilde{U})$ to any element of $\widetilde{R} := \Repirreven \sqcup \Repirreven \sqcup \Repirrodd \sqcup \Repirrodd$. We denote the $i$-th copy of $\Repirreven$ ($\Repirrodd$) by $\textup{Rep}_\textup{even}^{\textup{irr}, (i)}$ ($\textup{Rep}_\textup{odd}^{\textup{irr}, (i)}$). Let $V$ be a irreducible corepresentation in $\Repirreven$. Then $V$ and $z^* \cdot V \cdot z$ are corepresentations of $\widetilde{A}$. Actually, if $V$ is an irreducible subrepresentation of $U^{\boxtimes 2k}$ then $V$ is an irreducible subrepresentation of $(\widetilde{U} \boxtimes \overline{\widetilde{U}})^{\boxtimes k}$ and $z^* \cdot V \cdot z$ is an irreducible subrepresentation of $(\overline{\widetilde{U}} \boxtimes \widetilde{U})^{\boxtimes k}$. We consider $V$ as an element of $\Repirrevenone$ and $z^* \cdot V \cdot z$ as an element of $\Repirreventwo$. Similarly we see that if $V \in \Repirrodd$ then we can associate with it corepresentations $V \cdot z \in \Repirroddone$ and $z^* \cdot V \in \Repirroddtwo$. Note that elements $s$ from $\widetilde{S}$ give corepresentations $\widetilde{U}_s$ by this identification.
Consider a word $w = w_1 \dotsc w_k$  with letters in $\widetilde{R}$. We say that $w$ is reduced if in the sequence $\widetilde{U}_{w_1}, \dotsc, \widetilde{U}_{w_n}$ a $z$ is never followed by $z^*$ and $U_x$ is always followed by $z$ or $z^*$. In formal terms:
\begin{align*}
  \forall 1 \leq i \leq k - 1: \: & (w_i \in \Repirrevenone \cup \Repirroddtwo \Rightarrow w_{i + 1} \in \Repirreventwo \cup \Repirroddtwo) \wedge \\
    & (w_i \in \Repirreventwo \cup \Repirroddone \Rightarrow w_{i + 1} \in \Repirrevenone \cup \Repirroddone)
\end{align*}
Any such reduced word $w = w_1 \dotsc w_k$ gives rise to an irreducible corepresentation of $\widetilde{A}$ by $\widetilde{U}_w := \widetilde{U}_{w_1} \boxtimes \dotsc \boxtimes \widetilde{U}_{w_k}$ and different reduced words give rise to inequivalent corepresentations by \ref{theorem_irr_rep_of_free_products}. Since any iterated tensor product of $\widetilde{U}$ and $\overline{\widetilde{U}}$ decomposes as a sum of irreducible corepresentations of the type $\widetilde{U}_w$, where $w$ is a reduced word with letters in $\widetilde{R}$, any irreducible corepresentation of $\widetilde{A}$ is equivalent to some $\widetilde{U}_w$.

\begin{definition}
Consider now a word $w = w_1 \dotsc w_k$ with letters in $\widetilde{S}$. It is called connected if every $z$ is followed by a $z^*$. Formally:
\begin{align*}
  \forall 1 \leq i \leq k-1: \:	& (w_i \in \Sevenone \cup \Soddtwo \Rightarrow w_{i + 1} \in \Sevenone \cup \Soddone) \wedge \\ 
    & (w_i \in \Seventwo \cup \Soddone \Rightarrow w_{i + 1} \in \Seventwo \cup \Soddtwo)
\end{align*}
\end{definition}

The following definition says how we can associate irreducible corepresentations of $\widetilde{A}$ to words with letters in $\widetilde{S}$.

\begin{definition}
  If $w$ is an arbitrary word with letters in $\widetilde{S}$ then it has a unique decomposition $w = x_1 \dotsc x_l$ into maximal connected words. This gives rise to a unique reduced word $w'$ with letters in $\widetilde{R}$. We set $\widetilde{U}_w := \widetilde{U}_{w'}$
\end{definition}

Next we have to do some preparations in order to prove theorem  \ref{theorem_fusion_rules_of_complexifications}.

\begin{definition}
  Let $x=x_1 \dotsc x_m$ be a word in $\widetilde{S}$. Then $\check{x_i}$ is the letter in $S$ corresponding to $x_i$ and $\check{x} := \check{x_1}\check{x_2} \dotsc \check{x_m}$.
\end{definition}

\begin{remark}\label{remark_irreps_from_connected_words}
  Note that if $x$ is a connected word with letters in $S$ then according to remark \ref{remark_elements_of_S_tilde} it can be written as $z^{i_0} \cdot \check{x} \cdot z^{i_1}$, $i_0,i_1 \in \{0,1,-1\}$ and we have $\widetilde{U}_x = z^{i_0} \boxtimes U_{\check{x}} \boxtimes z^{i_1}$. 
\end{remark}

\begin{definition}
Let $x,y$ be connected words with letters in $\widetilde{S}$. We say that $(x,y)$ fits together if $xy$ is a connected word.
\end{definition}

\begin{lemma}\label{lemma_fusion_of_connected_words}
  Let $x = x_1 \dotsc x_m$ and $y=y_1 \dotsc y_n$ be connected words with letters in $\widetilde{S}$ such that $(x_m,y_1)$ fits together. Write $\widetilde{U}_x = z^{i_0} \boxtimes U_{\check{x}} \boxtimes z^{i_1}$ and $\widetilde{U}_y = z^{j_0} \boxtimes U_{\check{y}} \boxtimes z^{j_1}$. Then 
  \[\widetilde{U}_x \boxtimes \widetilde{U}_y = z^{i_0} \boxtimes \left (\sum_{x = ac, y = \overline{c}b}{U_{\check{a}\check{b}} \boxplus U_{\check{a} \cdot \check{b}}} \right) \boxtimes z^{j_1} = \sum_{x = ac, y = \overline{c}b}{\widetilde{U}_{ab} \boxplus \widetilde{U}_{a \cdot b}} .\]
\end{lemma}
\begin{proof}
  Since $(x,y)$ fits together, we have $z^{i_1} \boxtimes z^{j_0} = 1$. So by remark \ref{remark_fusion_rules_in_free_products} the first equation follows. We have to prove that for all $x = ac$, $y= \overline{c}b$
  \begin{enumerate}
    \item $z^{i_0} \boxtimes U_{\check{a}\check{b}} \boxtimes z^{j_1} = \widetilde{U}_{ab}$
    \item $z^{i_0} \boxtimes U_{\check{a} \cdot \check{b}} \boxtimes z^{j_1} = \widetilde{U}_{a\cdot b}$.
  \end{enumerate}
  In order to prove (1), note that $ab$ is connected, since $a,b$ are connected and $(a,b)$ fits together. So (1) follows from the way irreducible corepresentations are associated to connected words remarked in \ref{remark_irreps_from_connected_words}. \\
  For (2) note that, since $(a,b)$ fits together, $\check{a}\cdot \check{b} = \emptyset$ if and only if $a\cdot b = \emptyset$. If $a\cdot b \neq \emptyset$ then it is connected and (2) follows by remark \ref{remark_irreps_from_connected_words} again.
\end{proof}

Now we can give the proof of Theorem \ref{theorem_fusion_rules_of_complexifications}
\begin{proof}[Proof of Theorem \ref{theorem_fusion_rules_of_complexifications}]
  Let $x = x_1 \dotsc x_k$ and $y=y_1 \dotsc y_l$ be words with letters in $\widetilde{S}$. We have to show that
  \[\widetilde{U}_x \boxtimes \widetilde{U}_y = \sum_{x = ac, y=\overline{c}b}{\widetilde{U}_{ab} \boxplus \widetilde{U}_{a \cdot b}}\]
  Let $x=u_1 \dotsc u_m$ and $y=v_1 \dotsc v_n$ be the decomposition in maximal connected words. We identify them with letters in $\widetilde{R}$.  Then
  \begin{align*}
    \widetilde{U}_x	& = z^{i_0} \boxtimes U_{\check{u_1}} \boxtimes z^{i_1} \boxtimes U_{\check{u_2}} \boxtimes z^{i_2} \boxtimes \dotsm \boxtimes U_{\check{u_{m-1}}} \boxtimes z^{i_{m-1}} \boxtimes \underbrace{z^{i_m} \boxtimes U_{\check{u_m}} \boxtimes z^{i_{m+1}}}_{=\widetilde{U}_{u_m}}, \\
    \widetilde{U}_y	& = \underbrace{z^{j_0} \boxtimes U_{\check{v_1}} \boxtimes z^{j_1}}_{= \widetilde{U}_{v_1}} \boxtimes z^{j_2} \boxtimes U_{\check{v_2}} \boxtimes z^{j_3} \boxtimes \dotsm \boxtimes z^{j_{n-1}} \boxtimes U_{\check{v_{n-1}}} \boxtimes z^{j_n} \boxtimes U_{\check{v_n}} \boxtimes z^{j_{n+1}}
  \end{align*}
  with $i_1,...,i_{m-2},j_3,...,j_n \in \{1,*\}$, $i_0, i_m, j_0, j_2 \in \{0,*\}$ and $i_{m-1}, i_{m+1},j_1, j_{n+1} \in \{0, 1\}$. \\
  We are going to consider the two cases $(x_k,y_1)$ do or do not fit together. Assume that $(x_k, y_1)$ do not fit together. This means $z^{i_{m+1}} \cdot z^{j_0} \neq 1$. Then $\widetilde{U}_x \boxtimes \widetilde{U}_y$ is irreducible by Theorem \ref{theorem_irr_rep_of_free_products}. Moreover, $xy = u_1 \dotsc u_mv_1 \dotsc v_n$ is a decomposition in maximal connected words. So $\widetilde{U}_x \boxtimes \widetilde{U}_y = \widetilde{U}_{xy}$. On the other hand $(x_k, y_1)$ not fitting together implies $x_k \neq \overline{y_1}$ and $x_k \cdot y_1 = \emptyset$. So $\sum_{x = ac, y=\overline{c}b}{\widetilde{U}_{ab} \boxplus \widetilde{U}_{a \cdot b}} = \widetilde{U}_{xy}$. This completes the proof for the first case. \\
  Assume now that $(x_k,y_1)$ fits together. This means $z^{i_{m+1}} \cdot z^{j_0} = 1$. By Lemma \ref{lemma_fusion_of_connected_words}
  \begin{align*}
  \begin{split}
    \widetilde{U}_x \boxtimes \widetilde{U}_y & = z^{i_0} \boxtimes U_{\check{u_1}} \boxtimes z^{i_1} \boxtimes \dotsm \boxtimes U_{\check{u_{m-1}}} \boxtimes z^{i_{m-1}} \boxtimes \\
    & (\sum_{u_m = ac, v_1 = \overline{c}b}{\widetilde{U}_{ab} \boxplus \widetilde{U}_{a \cdot b}}) \boxtimes z^{j_2} \boxtimes U_{\check{v_2}} \boxtimes z^{j_3} \boxtimes \dotsm \boxtimes U_{\check{v_n}} \boxtimes z^{j_{n+1}} \\
    & = z^{i_0} \boxtimes U_{\check{u_1}} \boxtimes z^{i_1} \boxtimes \dotsm \boxtimes U_{\check{u_{m-1}}} \boxtimes z^{i_{m-1}} \boxtimes \\
    & ((\sum_{u_m = ac, v_1 = \overline{c}b, |a| \geq 1 \text{ or } |b| \geq 1}{\widetilde{U}_{ab} \boxplus \widetilde{U}_{a \cdot b}}) \boxplus \delta_{u_m,\overline{v_1}}\cdot 1) \boxtimes \\
    & z^{j_2} \boxtimes U_{\check{v_2}} \boxtimes z^{j_3} \boxtimes \dotsm \boxtimes U_{\check{v_n}} \boxtimes z^{j_{n+1}}.
  \end{split}
  \end{align*}
  By applying the induction hypothesis to the term
  \begin{gather*}
  z^{i_0} \boxtimes U_{\check{u_1}} \boxtimes \dotsm \boxtimes z^{i_{m-1}} \boxtimes \delta_{u_m,\overline{v_1}}\cdot 1 \boxtimes z^{j_2} \boxtimes U_{\check{v_2}} \boxtimes \dotsm \boxtimes z^{j_{n+1}} \\
  = \delta_{u_m,\overline{v_1}} \cdot \widetilde{U}_{u_1u_2 \dotsc u_{m-1}} \boxtimes \widetilde{U}_{v_2v_3 \dotsc v_n}
  \end{gather*}
  we obtain
  \[\widetilde{U}_x \boxtimes \widetilde{U}_y = \sum_{x = ac, y=\overline{c}b}{\widetilde{U}_{ab} \boxplus \widetilde{U}_{a \cdot b}}.\]
\end{proof}

We are now going to deduce the fusion rules of $\Akn$. The following result is proven in \cite{banicavergnioux09} and describes the fusion rules of $\Ahn$. 

\begin{theorem}
  Let $S_\textup{h} := \{u, p\}$ with fusion $u \cdot u = p \cdot p = p$, $u \cdot p = p \cdot u = u$ and trivial conjugation. The fusion rules of $(\Ahn, U_\textup{h})$ are given by the free fusion ring over $S_\textup{h}$ in such a way that $U_u \simeq U_\textup{h}$ and $U_p \boxplus 1 \simeq (u_{ij}^2)$.
\end{theorem}

Using this theorem we obtain the following corollary in the case $A = \Akn$.

\begin{corollary}
The irreducible corepresentations of $\Akn$ are described by the fusion set $S_{\textup{k}} := \{u,v,p,q\}$ with fusion given by
\begin{center}
\begin{tabular}{r||ccccc}
  $\cdot$	& $u$		& $v$		& $p$		& $q$		\\
  \hline \hline
  $u$		& $\emptyset$	& $q$		& $u$		& $\emptyset$	\\
  $v$		& $p$		& $\emptyset$	& $\emptyset$	& $v$		\\
  $p$		& $\emptyset$	& $v$		& $p$		& $\emptyset$	\\
  $q$		& $u$		& $\emptyset$	& $\emptyset$	& $q$		\\
\end{tabular}
\end{center}
and conjugation $\overline{u} = v$, $\overline{p} = p$, $\overline{q} = q$. \\
The elements of $S_{\textup{k}}$ correspond to the following corepresentations.
\begin{itemize}
  \item The class of the fundamental corepresentation $U$ is $U_u$.
  \item The class of $\overline{U}$ is $U_v$.
  \item The class of the corepresentation $(u_{ij}^* \cdot u_{ij})$ is $U_p \boxplus 1$
  \item The class of the corepresentation $(u_{ij} \cdot u_{ij}^*)$ is $U_q \boxplus 1$
\end{itemize}
\begin{proof}
  We only have to prove the part about the concrete description of $U_u$, $U_v$, $U_p$ and $U_q$. The fact that $U_u$ is the class of the fundamental corepresentation is obvious from the construction. $U_v \simeq \overline{U}$ follows directly. \\
  It is easy to check that $(u_{ij}^* \cdot u_{ij})$ and $(u_{ij} \cdot u_{ij}^*)$ are corepresentation of $\Akn$. We have the decomposition $U \boxtimes \overline{U} \simeq U_{uv} \boxplus U_p \boxplus 1$. Moreover the construction in this section shows that $U_{uv}$ is $n^2 - n$ dimensional and $U_p$ is $n-1$ dimensional. Since $(u_{ij} \cdot u_{ij}^*)$ is non trivial, it suffices to give at least two linearly independent intertwiners from the $n$ dimensional corepresentation $(u_{ij} \cdot u_{ij}^*)$ to $U \boxtimes \overline{U}$. Two such intertwiners are $\CC^n \rightarrow (\CC^n)^{\otimes 2}: \: e_i \mapsto e_i \otimes e_i$ and $\CC^n \rightarrow (\CC^n)^{\otimes 2}: \: e_i \mapsto \sum_j{e_j \otimes e_j}$. \\
  The proof for $(u_{ij} \cdot u_{ij}^*)$ works similarly.
\end{proof}
\end{corollary}

\bibliographystyle{plain}

\begin{thebibliography}{1}


\bibitem{banica96}
T. Banica,
\newblock {The representation theory of the free $O(n)$ compact quantum group.
  (Th\'eorie des repr\'esentations du groupe quantique compact libre $O(n)$.)}.
\newblock {\em C. R. Acad. Sci.}, (3) {\bf 322} (1996), 241-244.

\bibitem{banica97}
T. Banica,
\newblock {The free compact quantum group $U(n)$. (Le groupe quantique compact
  libre $U(n)$.)}.
\newblock {\em Commun. Math. Phys.}, (1) {\bf 190} (1997), 143-172.

\bibitem{banica99}
T. Banica,
\newblock {Symmetries of a generic coaction.}
\newblock {\em Math. Ann.}, (4) {\bf 314} (1999), 763-780.

\bibitem{banica08}
T. Banica,
\newblock {A note on free quantum groups. (Une note sur les groupes quantiques
  libres.)}.
\newblock {\em Ann. Math. Blaise Pascal}, (2) {\bf 15} (2008), 135-146.

\bibitem{banicaspeicher09}
T. Banica and R. Speicher,
\newblock {Liberation of orthogonal Lie groups.}
\newblock {\em Adv. Math.}, (4) {\bf 222} (2009), 1461-1501.

\bibitem{banicavergnioux09}
T. Banica and R. Vergnioux,
\newblock {Fusion rules for quantum reflection groups.}
\newblock {\em J. Noncommut. Geom.}, (3) {\bf 3} (2009), 327-359.

\bibitem{wang95}
S. Wang,
\newblock {Free products of compact quantum groups.}
\newblock {\em Commun. Math. Phys.}, (3) {\bf 167} (1995), 671-692.

\bibitem{wang98}
S. Wang,
\newblock {Quantum symmetry groups of finite spaces.}
\newblock {\em Commun. Math. Phys.}, (1) {\bf 195} (1998), 195-211.

\bibitem{woronowicz91}
S.L. Woronowicz,
\newblock {A remark on compact matrix quantum groups.}
\newblock {\em Lett. Math. Phys.}, (1) {\bf 21} (1991), 35-39.

\end{thebibliography}

% \providecommand{\bysame}{\leavevmode\hbox to3em{\hrulefill}\thinspace}
% \providecommand{\MR}{\relax\ifhmode\unskip\space\fi MR }
% % \MRhref is called by the amsart/book/proc definition of \MR.
% \providecommand{\MRhref}[2]{%
%   \href{http://www.ams.org/mathscinet-getitem?mr=#1}{#2}
% }
% \providecommand{\href}[2]{#2}

\end{document}